\documentclass[10pt]{amsart}

\usepackage{amssymb}
\usepackage[all]{xy}
\usepackage{mathrsfs}
\usepackage{enumerate}
\usepackage{tikz-cd}

\usepackage{bbm}

\usepackage{hyperref}
\hypersetup{
    colorlinks=true,
    linkcolor=blue,
    filecolor=magenta,      
    urlcolor=cyan,
}
\hypersetup{linktocpage}

\newcommand{\largewedge}{\mbox{\Large $\wedge$}}
\makeatletter
\newcommand*{\doublerightarrow}[2]{\mathrel{
  \settowidth{\@tempdima}{$\scriptstyle#1$}
  \settowidth{\@tempdimb}{$\scriptstyle#2$}
  \ifdim\@tempdimb>\@tempdima \@tempdima=\@tempdimb\fi
  \mathop{\vcenter{
    \offinterlineskip\ialign{\hbox to\dimexpr\@tempdima+1em{##}\cr
    \rightarrowfill\cr\noalign{\kern.5ex}
    \rightarrowfill\cr}}}\limits^{\!#1}_{\!#2}}}
\newcommand*{\triplerightarrow}[1]{\mathrel{
  \settowidth{\@tempdima}{$\scriptstyle#1$}
  \mathop{\vcenter{
    \offinterlineskip\ialign{\hbox to\dimexpr\@tempdima+1em{##}\cr
    \rightarrowfill\cr\noalign{\kern.5ex}
    \rightarrowfill\cr\noalign{\kern.5ex}
    \rightarrowfill\cr}}}\limits^{\!#1}}}
\makeatother

\newcommand{\nc}{\newcommand}

\nc\wh{\widehat}

\nc\on{\operatorname}

\nc\Gr{\on{Gr}}

\nc\Fl{\on{Fl}}

\newcommand{\limto}{{\displaystyle\lim_{\longrightarrow}}}
\newcommand{\rightlim}{\mathop{\limto}}

\newcommand{\leftlim}{\mathop{\displaystyle\lim_{\longleftarrow}}}
\newcommand{\limfromn}{\leftlim\limits_{\raise3pt\hbox{$n$}}}
\newcommand{\limton}{\rightlim\limits_{\raise3pt\hbox{$n$}}}

\newcommand{\rightlimit}[1]{\mathop{\lim\limits_{\longrightarrow}}\limits%
                    _{\raise3pt\hbox{$\scriptstyle #1$}}}

\newcommand{\leftlimit}[1]{\mathop{\lim\limits_{\longleftarrow}}\limits%
                    _{\raise3pt\hbox{$\scriptstyle #1$}}}

\makeatletter

\newcommand{\Rmnum}[1]{\expandafter\@slowromancap\romannumeral #1@}
\makeatother

\theoremstyle{definition}

\numberwithin{equation}{section}

\begin{document}

\title{Chern classes via derived determinant}

\author{Gleb Terentiuk}

\address{National research university "Higher school of Economics", Russian Federation}
\email{gleb.terentiuk@gmail.com}

%and
%\affilnum{2}Complete Second Author Address}

% Address / e-mail address of corresponding author
%\correspdetails{corr.email@math.edu}

%\date{}

\begin{abstract}
Motivated by the Chern-Weil theory, we prove that for a given vector bundle $E$ on a smooth scheme $X$ over a field $k$ of any characteristic, the Chern classes of $E$ in the Hodge cohomology can be recovered from the Atiyah class. Although this problem was solved by Illusie in \cite{i}, we present another proof by means of derived algebraic geometry.

Also, for a scheme $X$ over a field $k$ of characteristic $p$ with a vector bundle $E$ we construct elements $c^{cris}_n (E, \alpha(E)) \in H_{dR}^{2n} (X) $ using an obstruction $\alpha(E)$ to a lifting of $F^* E$ to a crystal modulo $p^2$ and prove that $c^{cris}_n (E, \alpha(E)) = n! \cdot c_{n}^{dR} (E)$, where $c_{n}^{dR} (E)$ are the Chern classes of $E$ in the de Rham cohomology and $F$ is the Frobenius map.
\end{abstract}

\maketitle
\tableofcontents

\section{Introduction}
\subsection{Motivation}
Let $X$ be a smooth scheme over a field $k$ and let $E$ be a vector bundle on $X$. Grothendieck defined \cite{g} the Chern classes $$c_n (E) \in F^n H_{dR}^{2n} (X), $$
where $F^n H_{dR}^{2n} (X)$ is the hypercohomology of the truncated de Rham complex: $F^n H_{dR}^{2n} (X) := R^{2n} \Gamma (X, 0 \to \cdots \to 0 \to \Omega_{X}^{n} \to \cdots \to \Omega_{X}^{\operatorname{dim}X}).$ We denote by $c_{n}^{h} (E) \in H^n (X, \Omega_{X}^{n})$ the image of $c_n (E)$ under the natural map $F^n H_{dR}^{2n} (X) \to H^n (X, \Omega_{X}^n).$
\par
Let $\operatorname{Conn} (E)$ be the sheaf of connections on $E;$ this is naturally a torsor over $\Omega_{X}^{1} \otimes \operatorname{End}(E).$ Denote by $\operatorname{At} (E) \in H^1 (X, \Omega_{X}^1 \otimes \operatorname{End} (E))$ the corresponding cohomology class. The first question addressed in this paper is whether $c_{n}^{h} (E)$ can be recovered from $\operatorname{At}(E).$ When $k = \mathbb{C}$ the (positive) answer follows easily from Chern-Weil theory. For example, $$c_{1}^{h} (E) = \operatorname{tr} A, ~ c_{2}^{h} (E) = \dfrac{1}{2}( \operatorname{tr}^2 A - \operatorname{tr} A^2), ~\text{and} ~ c_{3}^{h} (E) = \dfrac{1}{6}( 2\operatorname{tr} A^{3} + \operatorname{tr^{3}} A - 3\operatorname{tr} A^2 \operatorname{tr} A), $$ where $A$ stands for $\operatorname{At} (E)$ and by $A^n$ we mean the $n$-th power of $A$ in the algebra $H^* (X,  \Omega^1_{X} \otimes \operatorname{End}(E)).$ The same formulas hold for any field of characteristic $0$. However, if $\operatorname{char} k $ is positive the above formulas do not even make sense. Nevertheless, the answer is positive and first two parts of the paper are devoted to this question. Although this problem was solved by Illusie in \cite{i}, we present another proof by means of derived algebraic geometry. 
\par 
\subsection{The Main Results}
The first part contains necessary definitions and main tools meanwhile the proof itself is given in the second part. The main construction is as follows: for any smooth scheme $X$ with a locally free sheaf $L$, a perfect complex $\mathcal{E}$, and a morphism $A \in \operatorname{Hom}(\mathcal{E}, L[1] \otimes \mathcal{E})$, we construct an invertible element of degree zero $c(L, \mathcal{E}, A) \in H^{*} (X, \largewedge^{*} L) [t],$ 
where $H^{*} (X, \largewedge^{*} L) = \oplus_{n} H^{n} (X, \largewedge^{n} L)$ is a graded ring in a natural way and $\operatorname{deg} t = -1.$
\newline
This construction is based on the derived determinant morphism. Namely, we take the map $\pi: \operatorname{\underline{Spec}} \operatorname{S^{\bullet}} L[1] \to X,$ where $\operatorname{S^{\bullet}} L[1]$ is the free commutative simplicial $\mathcal{O}_X$-algebra constructed by $L[1]$ via the Dold-Kan correspondence and $$\operatorname{\underline{Spec}}: \operatorname{\textbf{sAlg}}_{\mathcal{O}_X} \to \operatorname{\textbf{dSch}}$$ is the relative spectrum functor. Thereafter, we interpret $A$ as an endomorphism of $\pi^* \mathcal{E}$. We denote the corresponding endomorphism by $tA$ and apply the derived determinant morphism to $\operatorname{id} - tA$ to obtain $c(L, \mathcal{E}, A)$ which we call the {\it characteristic polynomial} of $A$. The main theorem we are going to prove is the following.
\newline
\textbf{Theorem.} Let $X$ be a smooth scheme with a perfect complex $\mathcal{E},$ then $$c_k (\Omega_{X}^1, \mathcal{E}, \operatorname{At}(\mathcal{E})) = c_{k}^h (\mathcal{E}).$$
All we have to prove is that these classes coincide with Chern classes in the Hodge cohomology for line bundles, which is true since $\operatorname{tr}\operatorname{At}{\mathcal{E}} = c_{1}^{h} (\mathcal{E}),$ and also that the Whitney sum formula holds, which is true by virtue of the derived determinant construction.
\par
The third part is devoted to the following result. Let $X$ be a scheme over $\mathbb{F}_p$ and let $\mathcal{E}$ be a perfect complex on $X$. Notice that $F^* \mathcal{E}$ has a natural integrable connection with zero $p$-curvature and therefore can be viewed as an object of $\operatorname{Cris}(X/\mathbb{F}_p).$ There is an element $\alpha(\mathcal{E}) \in \operatorname{Ext}_{\operatorname{Cris}(X/\mathbb{F}_p)}^2 (F^* \mathcal{E}, F^* \mathcal{E})$ which is an obstruction to a lifting of $F^* \mathcal{E}$ to a crystal modulo $p^2.$ Denote by $c_{k}^{dR} (\mathcal{E})$ the image of $c_{k} (\mathcal{E})$ under the natural map $F^k H_{dR}^{2k} (X) \to H_{dR}^{2k} (X).$ It turns out that classes $k! \cdot c_{k}^{dR} (\mathcal{E})$ can be expressed in terms of $\alpha(\mathcal{E}).$ Similarly to the previous case, having $\mathcal{E}$ and $A \in  \operatorname{Ext}_{\operatorname{Cris}(X/\mathbb{F}_p)}^2 (F^* \mathcal{E}, F^* \mathcal{E}),$ we construct elements $c^{cris}_k (\mathcal{E}, A) \in H_{dR}^{2k} (X).$ Namely, denote the PD-envelope of the diagonal $X \xrightarrow{\Delta} X^n$ by $D(X^n)$, then a collection $D(X^n)$ with canonical projections forms a simplicial scheme $D(X).$ Then we interpret $A$ as an endomorphism $tA$ of $\pi^* F^* \mathcal{E},$ where $\pi: D(X) \times \operatorname{\underline{Spec}} \operatorname{S^{\bullet}} {\mathbb{Z}} [2] \to D(X).$ Afterwards, we apply the derived determinant morphism to $\operatorname{id} - tA$ to obtain an invertible element $c^{cris} (\mathcal{E}, A) \in \oplus_{k} H_{dR}^{2k} (X) \cdot \dfrac{t^k}{k!}.$ The main result that we are going to prove is the following.
\newline
\textbf{Theorem.} Let $X$ be a smooth scheme over $\mathbb{F}_p$ and let $\mathcal{E}$ be a perfect complex on $X$, then $$c^{cris}_k (\mathcal{E}, \alpha(\mathcal{E})) = k! \cdot c_{k}^{dR} (\mathcal{E}). $$
\newline
Again, we have to prove that $c_{1}^{dR} (\mathcal{E}) = \alpha(\mathcal{E})$ for line bundles and check that the Whitney sum formula holds. Although the second condition is a corollary of the splitting principle and properties of the derived determinant morphism, to prove the first condition we  reformulate it in the terms of an equivalence of torsors of sheaves.

 {\bf Acknowledgments.} I am grateful to Vadim Vologodsky for suggesting the problem, explaining the main ideas and reviewing the paper. Also, I would like to thank Alexander Petrov for explaining some ideas at the beginning of the work and Ezra Getzler for the interest to the paper. 

\section{Preliminaries} 
\subsection{Atiyah Class} Let $X$ be a smooth scheme over a field $k$ and let $\mathcal{E}$ be a coherent sheaf on $X$. Consider the diagonal embedding $\Delta: X \to X \times X,$ the corresponding ideal sheaf $\mathcal{I}$ and two canonical projections $X \times X\doublerightarrow{p}{q} X.$ Denote $\mathcal{O}_\Delta := \mathcal{O}_{X \times X}/\mathcal{I}$ and $\mathcal{O}_{\Delta^{(1)}} := \mathcal{O}_{X \times X} / \mathcal{I}^2.$ There is a short exact sequence
\begin{center}
$0 \to \mathcal{I}/\mathcal{I}^2 \to \mathcal{O}_{\Delta^{(1)}} \to \mathcal{O}_{\Delta} \to 0,$ 
\end{center}
and after tensoring it with $q^* \mathcal{E}$ and applying $p_*$ we get $$0 \to \mathcal{E} \otimes \Omega^1_X \to p_*(q^* \mathcal{E} \otimes \mathcal{O}_{\Delta^{(1)}}) \to \mathcal{E} \to 0.$$ We call the corresponding class of this extension $\operatorname{At} \mathcal{E} \in \text{Ext}^1 (\mathcal{E}, \mathcal{E} \otimes \Omega^1_X)$ the Atiyah class. 
\newline
More generally, if $E^{\bullet}$ is a finite complex of locally free sheaves, then we can similarly define the Atiyah class $\operatorname{At} E^{\bullet} \in \operatorname{Hom}(E^{\bullet}, E^{\bullet}[1] \otimes \Omega^1_X).$ 

\subsection{Derived Symmetric Powers}
Let $R$ and $S$ be commutative rings and let $F$ be a right exact functor $F: R$-$\operatorname{mod} \to S$-$\operatorname{mod}.$ When one wants to formally define left derived functors $L_i F,$ one takes a module $M$ over $R,$ chooses a projective resolution $P^{\bullet}$ and computes homology of $F(P^{\bullet})$.
Although a projective resolution is not unique, it is unique up to chain-homotopy equivalence. Thus, $F$ should preserve chain-homotopy equivalence. Unfortunately, the symmetric power functor is nonlinear. Dold and Puppe overcame this problem and defined the derived functors for the symmetric power functor, for more details see \cite{dp}. The idea is to use the Dold-Kan correspondence which states that there is an equivalence between the category of nonnegatively graded chain complexes and the category of simplicial abelian groups. Furthermore, this correspondence takes chain homotopies to simplicial homotopies. Therefore, the right recipe should look in the following way: make a complex into a simplicial abelian group, apply the symmetric power functor term-wise, then make it into a complex and compute its homology.
\newline
Important results that we will use state that if $\mathcal{M}$ is a flat module, then $LS^p(\mathcal{M}[1]) = \largewedge^p \mathcal{M}[p]$ and $L\largewedge^p(\mathcal{M}[1]) = \Gamma^p \mathcal{M}[p],$ see [\cite{i}, 4.3.4, 4.3.5].
\subsection{Determinant Map} For any locally free sheaf of a finite rank $\mathcal{E}$ on $X$ define $\operatorname{det} \mathcal{E}.$ Namely, if $U$ is an affine open subset of $X,$ then $\left.\left(\operatorname{det} \mathcal{E}\right)\right|_{U} :=\left(\largewedge^{\text{rk}~ \mathcal{E}} \mathcal{E}(U)\right)^{\sim}.$ 
\newline
In fact, there is a morphism of derived stacks $\operatorname{det_{Perf}}: \textbf{Perf} \to \textbf{Pic},$ see \cite{stv}.
An important result that we will use states that if $E$ and $F$ are perfect complexes, then $\operatorname{det_{Perf}}(E \oplus F) = \operatorname{det_{Perf}}(E) \otimes \operatorname{det_{Perf}}(F).$

\section{Chern Classes in Hodge Cohomology} 
\subsection{Characteristic Polynomial}
Let $X$ be a smooth scheme with a locally free sheaf $L$. Given a perfect complex $\mathcal{E}$ and $A \in \operatorname{Hom}(\mathcal{E}, L[1] \otimes \mathcal{E})$, we construct an element $c(L, \mathcal{E}, A) \in \left(H^{*} (X, \largewedge^{*} L)[t] \right)^*,$ where $H^{*} (X, \largewedge^{*} L) = \oplus_{k} H^k (X, \largewedge^{k} L)$ is a graded ring in a natural way and $\operatorname{deg} t = -1.$
\newline 
Make a complex $L[1]$ into a simplicial $\mathcal{O}_X$-module and apply the symmetric power functor term-wise, denote this free commutative simplicial $\mathcal{O}_X$-algebra by $\operatorname{S^{\bullet}} {L}[1].$ Let $\operatorname{\underline{Spec}}: \operatorname{\textbf{sAlg}}_{\mathcal{O}_X} \to \operatorname{\textbf{dSch}}$ be the relative spectrum functor, then denote $Y : = \operatorname{\underline{Spec}} \operatorname{S^{\bullet}} {L} [1].$  We have the map $\pi: Y \to X.$ 
\newline
\hypertarget{l:1}{\textbf{Lemma 1.1.}} An element $A \in \operatorname{Hom}(\mathcal{E}, L[1] \otimes \mathcal{E})$ provides us with an endomorphism $tA$ of $\pi^* \mathcal{E}$.
\newline
\textit{Proof.} To construct such an endomorphism it suffices to construct a morphism $\mathcal{E} \to \pi_* \pi^* \mathcal{E}$ by the adjunction. Identify $\pi_* \pi^* \mathcal{E}  \cong \mathcal{E} \otimes \operatorname{S^{\bullet}} L[1]$ via the projection formula. Take a natural morphism $i: \mathcal{E} \otimes L[1] \to \mathcal{E} \otimes \operatorname{S^{\bullet}} L[1]$ and the desirable morphism is $i \circ A: \mathcal{E} \to \pi_* \pi^* \mathcal{E}.$ $\qed$
\newline
\textbf{Lemma 1.2.} A morphism $\operatorname{id} - tA$ is an automorphism.
\newline
\textit{Proof.} It suffices to check that $\operatorname{id} - tA$ is an automorphism on the truncation of $Y$, but $tA$ is the zero morphism in $\pi_0(\operatorname{S^{\bullet}} L[1])$ because it shifts the filtration on $\pi^* \mathcal{E}.$ $\qed$ 
\newline
Since $\operatorname{id} - tA$ is an automorphism, it is possible to apply the determinant map, i.e. $\operatorname{det_{Perf}}({\operatorname{id}-tA})$ is an automorphism of $\operatorname{det}{\pi^* \mathcal{E}}$ which is a line bundle, therefore $\operatorname{det_{Perf}}({\operatorname{id}-tA}) \in \mathcal{O}^*_Y (Y).$
\newline
\textbf{Lemma 1.3.} One has $\mathcal{O}^*_Y (Y) = \left(H^* (X, \largewedge^* L) [t] \right)^*.$
\newline
\textit{Proof.} Homotopy groups of $\operatorname{S^{\bullet}} L[1]$ are $\largewedge^p L[p],$ then $\mathcal{O}_Y (Y) = \mathbb{H}^0 (X, \oplus_{p} \largewedge^p L[p]) = {\oplus_{p} H^p (X, \largewedge^p L) \cdot t^p},$ therefore $$\mathcal{O}^*_Y (Y) = \left(\oplus_{k} H^k (X, \largewedge^k L) \cdot t^k \right)^* = H^0 (X, \mathcal{O}^*_{X}) \oplus H^1 (X, L) \cdot t \oplus H^2 (X, \largewedge^2 L) \cdot t^2 \oplus \dotsc. \qed$$ 
\newline
Thus, we obtain elements $c_i (L, \mathcal{E}, A) \in H^i (X, \largewedge^i L).$
  
 \subsection{The Main Result} Apply the construction [2.1] in case of $L = \Omega^1_X$ and $A = \operatorname{At}\mathcal{E}.$ We claim that $c_i(\Omega^1,\mathcal{E}, \operatorname{At}\mathcal{E}) \in H^i(X, \Omega^i_X)$ are the Chern classes. By virtue of the construction the claim is true for line bundles, since $\operatorname{At}\mathcal{E} = c_1 (\mathcal{E})$ for line bundles, therefore it is true in general because for any perfect complexes $\mathcal{E}_1, \mathcal{E}_2$ we have $$c(L, \mathcal{E}_1 \oplus \mathcal{E}_2, A) = c(L, \mathcal{E}_1, A) \cdot c(L, \mathcal{E}_2, A), $$ i.e. the Whitney sum formula holds. 
 
\section{Chern Classes and Lifts Modulo $p^2$}
Let $X$ be a smooth scheme over $\mathbb{F}_p$ and denote by $F$ the Frobenius map. Then we have the direct image and the inverse image functors between categories of $\mathbb{F}_p$- crystals and $\mathbb{Z}_p$- crystals: \begin{tikzcd}[every arrow/.append style={shift left}]
 \operatorname{Cris}(X/\mathbb{F}_p) \arrow{r}{i_{cris *}} &\operatorname{Cris}(X/\mathbb{Z}_p) \arrow{l}{{i^*_{cris}}} 
 \end{tikzcd}. Note that for any perfect complex $\mathcal{E}$ there is a canonical connection $\nabla_{can}$ on $F^* \mathcal{E}$. Namely, let $\mathcal{E}$ be any $\mathcal{O}_X$-module, then for a local section of $F^* \mathcal{E}$ of the form $f \otimes s$ where $f$ is a function and $s$ is a local section of $\mathcal{E}$ define $\nabla_{can} (f \otimes s) = s \otimes df.$ This definition extends to complexes because $F$ is a flat morphism, i.e. $F^*$ is an exact functor.
Thus, any perfect complex $\mathcal{E}$ has a canonical integrable connection on $F^* \mathcal{E}$ with zero $p$-curvature, therefore we may view $F^* \mathcal{E}$ as an object of $\operatorname{Cris}(X/\mathbb{F}_p).$ According \cite{v} for any perfect complex $\mathcal{E}$ there is an exact triangle
\begin{center}
    $F^* \mathcal{E}[1] \to  (Li^*_{cris} \circ i_{cris *} ) F^* \mathcal{E} \to {F^* \mathcal{E}} \xrightarrow{\alpha(\mathcal{E})} F^* \mathcal{E}[2]$, 
\end{center}
i.e. an element $\alpha(\mathcal{E}) \in \operatorname{Ext}_{\operatorname{Cris}(X/\mathbb{F}_p)}^2 (F^* \mathcal{E}, F^* \mathcal{E}).$ Moreover, $\alpha(\mathcal{E}) = 0$ if and only if $F^* \mathcal{E}$ admits a lift to a crystal modulo $p^2.$
\newline
Let $\mathcal{E}$ be a perfect complex. Make a complex $\mathbb{Z} [2]$ into a free commutative simplicial algebra  $\operatorname{S^{\bullet}} \mathbb{Z} [2]$ and apply the relative spectrum functor to obtain a derived scheme $Y : = \operatorname{\underline{Spec}} \operatorname{S^{\bullet}} {\mathbb{Z}} [2].$ Consider the projection $\pi: Z := D(X) \times Y \to D(X) $ and take $\pi^* F^* \mathcal{E}$ where we view $F^* \mathcal{E}$ as a perfect complex on $D(X)$. Similarly as in the $\hyperlink{l:1}{\text{Lemma 1.2}}$ we obtain an automorphism $\text{id} - t \alpha(\mathcal{E})$ of $\pi^* F^* \mathcal{E},$ then we obtain $\operatorname{det_{Perf}}(\text{id} - t \alpha(\mathcal{E})) \in \mathcal{O}_{Z}^* (Z).$
\newline
$\textbf{Lemma 2.1}.$ One has $\mathcal{O}_{Z}^* (Z) = \Big(\oplus_{n} H^{2n}_{dR} (X) \cdot \dfrac{t^n}{n!} \Big)^*.$
\newline
\textit{Proof.} Using that $LS^n M[2] = L\Gamma^n M [2n]$ we obtain $$\mathcal{O}_Z(Z)  =\mathbb{H}^0(\oplus_{k} ( \mathcal{O}_{D(X)} \otimes \Gamma^k \mathbb{Z}[2k])) = \oplus_{k} H^{2k} (D(X), \mathcal{O}_{D(X)}) \otimes \Gamma^k \mathbb{Z} = \oplus_{k} H^{2k}_{dR} (X) \cdot \dfrac{t^k}{k!}$$ and we obtain the desirable result. For the last isomorphism see \cite{bo}. $\qed$
\newline
Thus, we have $\operatorname{det_{Perf}} (\operatorname{id} - t\alpha(\mathcal{E})) = \sum c^{{cris}}_k (\mathcal{E}, \alpha(\mathcal{E})) \dfrac{t^k}{k!}$ where $c^{{cris}}_k (\mathcal{E}, \alpha(\mathcal{E})) \in H^{2k}_{dR} (X).$

Suppose that $\mathcal{E}$ is a line bundle. The sheaf of connections $\operatorname{Conn}(E)$ is naturally a sheaf of torsors over $\Omega^1_X$ and its subsheaf of flat connections $\operatorname{FConn}(E)$ is a sheaf of torsors over $\Omega^1_{closed}.$   
\newline 
\textbf{Lemma 2.2.} A bundle $F^* \mathcal{E}$ can be locally lifted modulo $p^2.$
\newline
\textit{Proof.} Take any local lift of $X$ modulo $p^2$, then any trivialization of $\mathcal{E}$ yields transition functions $f_\alpha.$ Transition functions of $F^* \mathcal{E}$ are $f^p_\alpha$ and they can be lifted modulo $p^2$ by taking the Teichmüller lifts.  $\qed$
\newline
Denote this lift by $(F^* \mathcal{E})_{can}.$ Note that this construction does not depend on the lift of $X$ in the following sense: if $X_1$ and $X_2$ are local lifts of $X$ modulo $p^2,$ $(F^* \mathcal{E})^1_{can}$ and $(F^* \mathcal{E})^2_{can}$ are the corresponding lifts of $F^* \mathcal{E},$ and if $f: X_1 \to X_2$ is any isomorphism such that it is the identity modulo $p,$ then $f^* (F^* \mathcal{E})^2_{can} = (F^* \mathcal{E})^1_{can}.$ 

Consider the sheaf $\operatorname{Lift}(\mathcal{E})$ of lifts $M$ of $F^* \mathcal{E}$ to a crystal modulo $p^2$ such that $M$ is isomorphic to $(F^* \mathcal{E})_{can}$ with a flat connection $J$ which coincides with $\nabla_{can}$ modulo $p.$ Thus, $\operatorname{Lift}(\mathcal{E})$ is a sheaf of torsors over $\Omega^1_{closed}.$ Images of $\operatorname{FConn}(\mathcal{E})$ and $\operatorname{Lift}(\mathcal{E})$ under the map $H^1 (X, \Omega_{closed}^1) \to H^2_{dR} (X)$ are $c_{1}^{dR} (E)$ and $\alpha (E)$ respectively. 
\newline
\textbf{Lemma 2.3.} Sheaves of torsors $\operatorname{FConn}(\mathcal{E})$ and $\operatorname{Lift}(\mathcal{E})$ are equivalent.
\newline
\textit{Proof.} Use [\cite{pvv}, Lemma 2.5]. $\qed$
\newline
Thus, we obtain $c^{{cris}}_n (\mathcal{E}, \alpha({\mathcal{E})}) = n! \cdot c^{\text{dR}}_n (\mathcal{E})$ from the splitting principle.

\end{document}